\documentclass[12pt]{article}

\usepackage{amsmath, amsthm, amssymb, enumerate}
\usepackage{color}
\usepackage{mathtools}  
\usepackage{datetime}
\usepackage{fancyhdr}
\usepackage{eufrak}
\chardef\coloryes=1 
\chardef\isitdraft=1 

\ifnum\isitdraft=1 
   \def\version{9} 
   \def\eqref#1{({\ref{#1}})}                

\definecolor{labelkey}{gray}{.3}

\definecolor{refkey}{rgb}{.3,0.3,0.3}

\else

  \def\startnewsection#1#2{\section{#1}\label{#2}\setcounter{equation}{0}}   
  \def\nnewpage{} 
\fi

 \usepackage{amsmath}
\usepackage{amsfonts}
\usepackage{amssymb}
\def\sg {{\sf g}}

\newcommand{\cP}{{\mathfrak P}}

\newcommand{\bex}{\bbE^{\bbP_0}_{x_0}}
\newcommand{\bext}{\bbE^{\bbP_0}_{t_i}}
\newcommand{\bextt}{\bbE^{\bbP_0}_{t}}
\newcommand{\bextpi}{\bbE^{\bbP_0}_{x_i}}
\newcommand{\bextq}{\bbE^{Q}_{t_i}}

\newtheorem{theorem}{Theorem}[section]
\newtheorem{lemma}{Lemma}[section]

\newtheorem{remark}{Remark}[section]
\newtheorem{definition}{Definition}[section]

\newtheorem{cor}{Corollary}[section]

\newcommand{\bgeqn}{\begin{eqnarray}}
\newcommand{\edeqn}{\end{eqnarray}}
\newcommand{\bgeq}{\begin{eqnarray*}}
\newcommand{\edeq}{\end{eqnarray*}}
\newcommand{\bec}{\begin{center}}
\newcommand{\enc}{\end{center}}
\newcommand{\beqn}{\begin{equation}}
\newcommand{\eeqn}{\end{equation}}
\newcommand{\belign}{\begin{align}}
\newcommand{\elign}{\end{align}}
\newcommand{\belem}{\begin{lemma}}
\newcommand{\elem}{\end{lemma}}
\newcommand{\bethm}{\begin{theorem}} 
\newcommand{\ethm}{\end{theorem}}
\newcommand{\beitem}{\begin{itemize}}
\newcommand{\eitem}{\end{itemize}}
\newcommand{\benum}{\begin{enumerate}}
\newcommand{\enum}{\end{enumerate}}
\newcommand{\beproof}{\begin{proof}}
\newcommand{\eproof}{\end{proof}}
\newcommand{\bedef}{\begin{definition}}
\newcommand{\edefn}{\end{definition}}

\newcommand{\trm}{\textrm}
\newcommand{\berem}{\begin{remark}}
\newcommand{\erem}{\end{remark}}
\newcommand{\pis}{ \pi^2_s (\sg_s^S)^2}

\newcommand{\bbR}{\mathbb{R}}
\newcommand{\bea}{\begin{eqnarray}}
\newcommand{\eea}{\end{eqnarray}}
\newcommand{\bearr}{\begin{eqnarray}}
\newcommand{\enarr}{\end{eqnarray}}

\newcommand{\Q}{{\cal Q}}

\newcommand{\F}{{\cal F}}

\newcommand{\C}{{\cal C}}

\newcommand{\B}{{\cal B}}

\newcommand{\bbP}{\mathbb{P}}
\newcommand{\bbN}{\mathbb{N}_0}

\newcommand{\be}{\begin{equation}}
\newcommand{\ee}{\end{equation}}
\newcommand{\rar}{\rightarrow}

\newcommand{\supi}{\sup_{\pi \in \Pi_{\mathbf{ad}}}}

\newcommand{\infir}{\inf_{Q \in \Q^c_{\ti}}}

\newcommand{\infit}{\inf_{Q \in \Q_{\ti}}}
\newcommand{\supit}{\sup_{\pi \in \Pi_{\mathrm{ad}}}}

\def\dist{\mathop{\rm dist}}
\def\w{\omega}

\def\O{\Omega}

\def\supp {{\rm supp}}

\def\bbP{{\mathbb{P}}}
\def\bbE{{\mathbb{E}}}

\def \rhot1{\rho_{t+1}}
\def \rhot{\rho_{t}}

\def \rhp{\rightharpoonup}
\def \sg {\sigma}
\def \ti {{[t_i,T)}}

\def \piti {T}
\def \tsg {\theta^\sigma}
\def \esg{\eta^\sigma}
\def \intzt {\int_0^t}
\def \intst {\int_s^t}

\def \tmu   {\theta^\mu}
\def \emu   {\eta^\mu}
\def \dpi { (\pi_s - \pi^n_s) }

\def \nut {\nu_t}
\def \gam {\gamma}

\def\beal#1\eal{\begin{align}#1\end{align}}

 \renewcommand{\theequation}{\thesection.\arabic{equation}}

\begin{document}
\def\ques{{\cor \underline{??????}\cob}}
\def\nto#1{{\coC \footnote{\em \coC #1}}}
\def\fractext#1#2{{#1}/{#2}}
\def\fracsm#1#2{{\textstyle{\frac{#1}{#2}}}}   
\def\nnonumber{}


\def\cor{{}}
\def\cog{{}}
\def\cob{{}}
\def\coe{{}}
\def\coA{{}}
\def\coB{{}}
\def\coC{{}}
\def\coD{{}}
\def\coE{{}}
\def\coF{{}}
\ifnum\coloryes=1

  \definecolor{coloraaaa}{rgb}{0.1,0.2,0.8}
  \definecolor{colorbbbb}{rgb}{0.1,0.7,0.1}
  \definecolor{colorcccc}{rgb}{0.8,0.3,0.9}
  \definecolor{colordddd}{rgb}{0.0,.5,0.0}
  \definecolor{coloreeee}{rgb}{0.8,0.3,0.9}
  \definecolor{colorffff}{rgb}{0.8,0.3,0.9}
  \definecolor{colorgggg}{rgb}{0.5,0.0,0.4}

 \def\cog{\color{colordddd}}
 \def\cob{\color{black}}
 \def\cor{\color{red}}
 \def\coe{\color{colorgggg}}

 \def\coA{\color{coloraaaa}}
 \def\coB{\color{colorbbbb}}
 \def\coC{\color{colorcccc}}
 \def\coD{\color{colordddd}}
 \def\coE{\color{coloreeee}}
 \def\coF{\color{colorffff}}

\fi
\ifnum\isitdraft=1
   \chardef\coloryes=1 
   \baselineskip=17pt
\pagestyle{myheadings}
\reversemarginpar

\def\const{\mathop{\rm const}\nolimits}  
\def\diam{\mathop{\rm diam}\nolimits}    

 \def\llabel#1{\label{#1}{\ \mbox{\rm\color{red} {#1}\color{black}}}}

\def\rref#1{{\ref{#1}{\rm \tiny \fbox{\tiny #1}}}}
\def\theequation{\fbox{\bf \thesection.\arabic{equation}}}
\def\ccite#1{{\cite{#1}{\rm \tiny ({#1})}}}

\def\startnewsection#1#2{\newpage\cog \section{#1}\cob\label{#2}

\setcounter{equation}{0}
\pagestyle{fancy}

\lhead{\cob Section~\ref{#2}, #1 }
\cfoot{}
\rfoot{\thepage}
\lfoot{\cob{\today,~\currenttime}~(c75-iklt2, Version~\fbox{\version})}}
\chead{}
\rhead{\thepage}
\def\nnewpage{\newpage}

\newcounter{startcurrpage}
\newcounter{currpage}

\def\llll#1{{\rm\tiny\fbox{#1}}}
   \def\blackdot{{\color{red}{\hskip-.0truecm\rule[-1mm]{4mm}{4mm}\hskip.2truecm}}\hskip-.3truecm}
   \def\bdot{{\coC {\hskip-.0truecm\rule[-1mm]{4mm}{4mm}\hskip.2truecm}}\hskip-.3truecm}
   \def\purpledot{{\coA{\rule[0mm]{4mm}{4mm}}\cob}}
   \def\pdot{\purpledot}
\else
   \baselineskip=15pt
   \def\blackdot{{\rule[-3mm]{8mm}{8mm}}}
   \def\purpledot{{\rule[-3mm]{8mm}{8mm}}}
   \def\pdot{}
\fi

\def\nts#1{{\hbox{\bf ~#1~}}} 
\def\nts#1{{\cor\hbox{\bf ~#1~}}} 
\def\ntsf#1{\footnote{\hbox{\bf ~#1~}}} 
\def\ntsf#1{\footnote{\cor\hbox{\bf ~#1~}}} 
\def\bigline#1{~\\\hskip2truecm~~~~{#1}{#1}{#1}{#1}{#1}{#1}{#1}{#1}{#1}{#1}{#1}{#1}{#1}{#1}{#1}{#1}{#1}{#1}{#1}{#1}{#1}\\}
\def\biglineb{\bigline{$\downarrow\,$ $\downarrow\,$}}
\def\biglinem{\bigline{---}}
\def\biglinee{\bigline{$\uparrow\,$ $\uparrow\,$}}

\mathtoolsset{showonlyrefs=true}

\def\Oti{\Omega_{[t_i,T]}}
\def\Qti{\Q_{[t_i,T]}}
\def\Qtic{\Qti^c}
\def\ti{[t_i,T]}
\def\tiz{[0,T]}
\def\tilde{\widetilde}
\def\pft{\cP(\O|\F_{t_i})}
\newcommand{\supq}{\sup_{\pi \in \Pi_{\mathbf{ad}}}}
\newcommand{\infq}{\inf_{Q \in \Qti}}
\newcommand{\infqc}{\inf_{Q \in \Qtic}}

\newtheorem{Theorem}{Theorem}[section]
\newtheorem{Corollary}[Theorem]{Corollary}
\newtheorem{Proposition}[Theorem]{Proposition}
\newtheorem{Lemma}[Theorem]{Lemma}
\newtheorem{Remark}[Theorem]{Remark}
\newtheorem{Example}[Theorem]{Example}
\newtheorem{Assumption}[Theorem]{Assumption}
\newtheorem{Claim}[Theorem]{Claim}
\newtheorem{Question}[Theorem]{Question}
\def\theequation{\thesection.\arabic{equation}}
\def\endproof{\hfill$\Box$\\}
\def\square{\hfill$\Box$\\}
\def\comma{ {\rm ,\qquad{}} }            
\def\commaone{ {\rm ,\qquad{}} }         
\def\dist{\mathop{\rm dist}\nolimits}    
\def\sgn{\mathop{\rm sgn\,}\nolimits}    
\def\Tr{\mathop{\rm Tr}\nolimits}    
\def\div{\mathop{\rm div}\nolimits}    
\def\supp{\mathop{\rm supp}\nolimits}    
\def\divtwo{\mathop{{\rm div}_2\,}\nolimits}    
\def\re{\mathop{\rm {\mathbb R}e}\nolimits}    
\def\div{\mathop{\rm{Lip}}\nolimits}   
\def\indeq{\qquad{}}                     
\def\period{.}                           
\def\semicolon{\,;} 


\title{Portfolio Optimization with Nondominated Priors and Unbounded Parameters}
\author{Kerem U\u{g}urlu}
\maketitle

\date{}

\begin{center}
\end{center}

\medskip

\indent Department of Applied Mathematics, University of Washington, Seattle, WA 98195\\
\indent e-mail:keremu@uw.edu

\begin{abstract}
We consider classical Merton problem of terminal wealth maximization in finite horizon. We assume that the drift of the stock is following Ornstein-Uhlenbeck process and the volatility of it is following GARCH(1) process. In particular, both mean and volatility are unbounded. We assume that there is Knightian uncertainty on the parameters of both mean and volatility. We take that the investor has logarithmic utility function, and solve the corresponding utility maximization problem explicitly. To the best of our knowledge, this is the first work on utility maximization with unbounded mean and volatility in Knightian uncertainty under nondominated priors.  
\end{abstract}

\noindent\thanks{\em Mathematics Subject Classification\/}:
91B28;93E20\\
\noindent\thanks{\em Keywords:\/}
Merton Problem; Knightian uncertainty; Robust optimization

\section{Introduction}

Starting with the pioneering works of \cite{M, B, AS, RM,  BSM}, the underlying risky assets are modelled as Markovian diffusions, where there exists a fixed underlying reference probability measure $\bbP$ that is retrieved from historical data of the price movements. However, it is mostly agreed that it is impossible to precisely identify $\bbP$. Hence, as a result, model ambiguity, also called Knightian uncertainty, in utility maximization is inevitably taken into consideration. Namely, the investor is diffident about the odds, and takes a robust approach to the utility maximization problem, where she minimizes over the priors, corresponding to different scenarios, and then maximizes over the investment strategies. 

The literature on robust utility maximization in mathematical finance, (see e.g. \cite{DW,MQ,AS,LSW,CF,CE,RO}), mostly assumes that the set priors is dominated by a reference measure $\bbP$. Hence, it presumes a setting where volatility of risky assets  are perfectly known, but drifts are uncertain. Namely, these approaches assume the equivalence of priors. In particular, they assume the equivalence of probability measures $P$ with a dominating reference prior $\bbP_0$. 
In this direction to mention some of the related works, \cite{GS} proposed to weaken
the strong independence axiom (also called sure thing principle used previously by \cite{LS} and \cite{AA}) to justify (subjective) expected utility. Later, \cite{ADEH} introduced coherent risk measures in the spirit of the construction as in \cite{GS}. This theory of coherent risk measures has been generalized in several directions afterwards,  (see e.g. \cite{MMDR,MMRU,FS,R,RS,ES,SL1,SL2} among others).

By contrast, we are studying the case, where the set of priors, denoted by $\cP(\O)$, are nondominated. Hence, there exists no dominating reference prior $\bbP_0$. Some of the related works in this direction are as follows: \cite{HS} studied the case, where uncertainty in the volatility is due to an unobservable factor; \cite{N} studied a similar setting in discrete time and has shown the existence of optimal portfolios; \cite{NN2, MQ} works in a jump-diffusion context, with ambiguity on drift, volatility and jump intensity. \cite{DK} establishes a minimax result and the existence of a worst-case measure in a setup where prices have continuous paths and the utility function is bounded. \cite{LR} works in a diffusion context, where uncertainty is modelled by allowing drift and volatility to vary in two constant order intervals. The optimization using power utility of the from $U(x) = x^\gamma$ for $0 < \gamma <1$ is then performed via a robust control (G-Brownian motion) technique, which requires the uncertain volatility matrix is restricted to be of diagonal type. We refer the reader to \cite{P} for a detailed exposure on G-Brownian motion and its applications. In a more general setting \cite{MN} works in a continuous time setting, where the stock prices are allowed to be general discontinuous semi-martingales and strategies are required to be compact. 

Contrary to the current literature, we assume that the drift of the stock process is modelled by Ornstein-Uhlenbeck (OU)-process, and volatility of the stock process is modelled by GARCH(1)-process. In particular, both mean and volatility of the stock take unbounded values in $\bbR$. We assume that there is Knightian uncertainty on the parameters of both mean and volatility. We further take the investor has logarithmic preference but being diffident about the underlying dynamic parameters wants to reconsider the optimal parameters at prespecified  small time intervals. To the best of our knowledge, this is the first work on utility maximization with unbounded mean and volatility in Knightian uncertainty under nondominated priors.  

The rest of the paper is as follows. In Section 2, we describe the model dynamics of the problem and state our general main problem. In Section 3, we solve the investor's problem using logarithmic utility. Finally, in Section 4, we discuss our results and conclude the paper. 
\section{Model Dynamics and Investor's Value Function}

\subsection{Framework}

We let $\O = C_0[0,T]$ be the set of all continuous paths $(\w_t)_{t\geq0}$ starting at $\w_0 = 0$ taking values in $\bbR$ over the finite time horizon $[0,T]$ equipped with uniform norm 
\beqn
\lVert \w \rVert_\infty = \sup_{t \in [0,T]} |\w_t|.
\eeqn
We further take the corresponding induced metric on $\O$ in the usual way
\beqn 
d(f,g):= \lVert f-g\rVert_\infty,
\eeqn and take the topology generated by the open sets in this metric as the topology of uniform convergence on $[0,T]$. We further define $W_t$ as the coordinate functional, i.e. for $\w \in \O$ and $t \in [0,T]$, we let $W_t(\w) := \w_t$. We denote the corresponding Borel sigma algebra on $\O$ by $\F_t = \sigma\{W_s(\w); s\leq t\}$, and the standard Wiener measure by $\bbP_0$ on $\O$. Namely, $\bbP_0$ is the unique measure on $\O$, which satisfies the following properties:
\begin{itemize}
\item $\bbP_0 (\w \in \Omega: W_0(\w) = 0)  =1$,
\item For any $f \in \C_b^\infty(\bbR)$, $\C_b^\infty(\bbR)$ being infinitely many times continuously differentiable bounded functions, the stochastic process
\begin{align}
(t,\w) \rar f(W_t(\w)) - \frac{1}{2}\int_0^t  \Delta f(W_s(\w))ds
\end{align}
is an $(\F_t, \bbP_0)$-martingale, 
\end{itemize}

We take the Wiener measure $\bbP_0$ as our \textit{reference measure} and $\cP(\O)$ as the Polish space equipped with weak topology of probability measures on $\O$. We further denote by $\cP(\O|\F_t)$ for $0 \leq t \leq T$ the Polish space of regular conditional probability measures (r.c.p.m.'s) on $\O$ (see \cite{D}, Thm 5.1.9). 
\subsection{Model Dynamics}
We consider a market consisting of one risky asset, $S_t$, and one riskless asset $R_t$. We assume that the semimartingale $S_t$ representing the price of the risky asset satisfies the following dynamics
\begin{align} 
S_0  &= s_0,\; s_0 > 0\\
\label{eqn250}
d S_t &= S_t(\mu_t dt +\sigma^S_t d W^S_t),\; \bbP_0\; \mathrm{ a.s.}.
\end{align}
Here $W^S_t$ and $W^\sg_t$ are independent $(\O,\F_t)$ Brownian motions, and both the drift term $\mu_t$ and the volatility $\sigma^S_t$ are stochastic processes satisfying the following dynamics 
\begin{align}
\label{eqnmu}
d\mu_t &= \theta^\mu_t(\eta_t^\mu - \mu_t)dt + \sigma^\mu dW^\mu_t,\\
\label{eqnsg}
d (\sigma^S_t)^2  &= \theta^\sigma_t(\eta_t^\sg - (\sigma^S_t)^2)dt + \xi(\sigma^S_t)^2dW^\sigma_t.
\end{align}
Here, $\theta_t^\mu,\eta_t^\mu, \theta_t^\sigma,\eta^\sg_t$ are piecewise constant on $[t_i,t_{i+1})$ for $i=0,\ldots,N$, with $t_0 =0$ and $t_N = T$ and strictly positive, whereas $\sg^\mu,\xi > 0$ are positive constants. Namely, we assume that the drift term, $\mu_t$, is an Ornstein-Uhlenbeck (OU) process, whereas the volatility $(\sg^S_t)^2$ is a GARCH(1) process. It is easy to see that the explicit solutions of $\mu_t$ and $(\sg^S_t)^2$ are 
\begin{align}
(\sg^S_t)^2 &= (\sg^S_0)^2\intzt \tsg_s\esg_s\exp\{ -\intst (\tsg_r + \frac{1}{2}\xi^2)dr + \xi(W_t - W_s) \}ds\\
\mu_t &= \mu_0 \exp\{-\intzt \tmu_s ds\} + \intzt \exp\{ -\intst \tmu_r dr\} \sg^\mu dW_s
\end{align}
We further assume that the riskless asset's dynamics satisfies 
$d R_t = rRdt$ with $r >0$.

\subsection{Model Uncertainty}
We denote by $\gamma \triangleq (\theta^\mu,\theta^\sg, \eta^\mu, \eta^\sg)$ as the quadruple of piecewise constant positive functions and let $\gamma \in \Gamma$ be a collection of them. We assume that $\Gamma$ is uniformly bounded, i.e. $|\gamma(t)| \leq M$ for some constant $M$ for all $\gamma \in \Gamma$ and $0 \leq t \leq T$. 

\begin{Remark} The quadratic variation terms of $\mu_t$  and $\sigma_t$ can be deduced with certainty from the path $\w$, but the drift parameters of these processes can only be estimated from historical data. Hence, it is natural to introduce uncertainty only in the corresponding drift terms of $\mu_t$ and $\sg_t$, whereas the quadratic variation terms $\xi,\sigma^\mu$ are fixed known constants. 
\end{Remark}

We next give the following a priori estimates that are to be used in the rest of the paper. 

\belem
\label{lem21}
Let $(\sg^S_t)^2$ and $\mu_t$ be as defined in \eqref{eqnmu} and \eqref{eqnsg}, respectively. Then, we have 
\beal
&\bbE^{\bbP_0}[\int_0^T\max_{\gamma \in \Gamma} (\sg^S_t)^{2n} dt] < \infty \\
&\bbE^{\bbP_0}[\int_0^T\max_{\gamma \in \Gamma} (\mu_t)^n dt] < \infty
\eal
for any $n \in \bbN$. In particular, 
\beal
&\bbE^{\bbP_0}_t[\int_{t}^T \max_{\gamma \in \Gamma} (\sg^S_s)^{2n} ds] < \infty\; \bbP_0\trm{-a.s.}\\
&\bbE^{\bbP_0}_t[\int_t^T \max_{\gamma \in \Gamma}(\mu_s)^n ds  ] < \infty \; \bbP_0\trm{-a.s.}
\eal
for any $t \in [0,T]$. 
\elem

\beproof
For ease of notation, denoting $(\sg^S_t)^2 = \nu_t$ and using constant $C$ interchangeably below, we have 
\beal
\nut &\leq \intzt|\tsg_s\esg_s|\exp\big(\intst (\tsg_r + \frac{1}{2}\xi^2)dr\big) + \xi|W_t -W_s|ds
\eal 
Since $\tsg,\esg$ and $\xi$ are uniformly bounded, we have 
\beqn
\max_{\gam} |\nut|^n \leq \intzt C \exp\bigg( C(t-s) + C|W_t - W_s| \bigg) ds
\eeqn

Hence, for any $n \geq 1$, via Jensen's inequality we have
\beqn
\max_{\gam}|\nut|^n \leq \intzt C \exp\bigg( C(1 + |W_t - W_s|) \bigg)ds.
\eeqn
Thus, 
\beqn
\bbE^{\bbP_0}[\max_{\gam}|\nut|^n] \leq C \intzt \bbE [\exp(C(1+ |W_t -W_s|))]ds < \infty
\eeqn
and 
\beal
\bbE[\int_0^T \max_{\gamma}|\nu_t|^n] < \infty.
\eal
For $\mu_t$ and $n \geq 1$, we have the following inequalities modulo $\bbP_0$-a.s.
\beal
|\mu_t| &\leq C + C|\int_0^t e^{-\int_s^t\theta_r^\mu}dW_s|\\
|\mu_t|^{2n} &\leq C + C\bigg( \bigg| \int_0^t e^{-\int_s^t\theta_r^\mu}dW_s \bigg|\bigg)^{2n}\\
\bbE^{\bbP_0}[\max_{\gamma \in \Gamma} |\mu_t|^{2n}] &\leq C + C\int_0^T e^{-2n\int_s^t\theta_r^\mu dr}ds\\
&\leq C(1+T)\\
&< \infty
\eal
Since $x^n \leq x^{2n} + 1$ for $x > 0$ and any $n \geq 1$, we have 
\beal
\bbE^{\bbP_0}[\max_{\gamma \in \Gamma} |\mu_t|^n] &< \infty \\
\bbE^{\bbP_0}[\int_0^T (\max_{\gamma \in \Gamma} |\mu_t|^n) dt] &< \infty. 
\eal 
Hence, we conclude the the proof. 
\eproof

\subsection{Alternative Models}

At each fixed time $t_i \in [0,T]$ with $0=t_0 < t_1 <\ldots < t_N = T$ we consider the class of alternative models denoted by $\Q_{[t_i,T]}$. These are the set of r.c.p.m.'s on $\O$
that are induced by the process $S(\cdot)$ as in Equation \eqref{eqn250}. Namely, 
\begin{align} 
\label{eqn26}
\Qti = \big\{ &P \in \cP(\O|\F_{t_i})|\; P \textrm{ is the r.c.p.m. induced by } \\
&S \textrm{ satisfying } d S_u =  S_u( \mu_u du +  \sigma^S_u d W^S_u) du,\; \bbP_0\mathrm{-a.s.}\;\\ 
&\textrm{for } u \in [t_i,T] \textrm{ given } \{S_r(\w),0\leq r \leq t_i\}\\
&\textrm{with } (\theta_u^\mu,\theta_u^\sigma,\eta_u^\sg,\eta_u^\mu)_{0\leq u \leq T} \in \Gamma
\big\}.
\end{align}

We note here that by Theorem 11.2 in \cite{RW}, there exists a strong solution to the Equation \eqref{eqn250} on $(\O, \F_T, \bbP_0)$. Namely, denoting $C_0[0,T] = \O$, there exists an $\F_T$ measurable mapping $G:\O \rar \O$ such that $X(\cdot) \equiv G(x_{t_i},W(\cdot))$ solves Equation \eqref{eqn250} on $(\O,\F_T, \bbP_0)$, as in Definition 10.9 in \cite{RW}. We further note that there is a one-to-one correspondence between the set of r.c.p.m.'s $\Qti$ and the compact set $\Gamma$. Namely, $\gamma = (\theta^\mu,\theta^\sigma,\eta^\sg\eta^\mu) \in \Gamma $ uniquely defines $S_t$ on $[0,T]$ $\bbP_0$-a.s. 
We denote the r.c.p.m. induced by $S_T$ for a fixed $\gamma \in \Gamma$ and $A \in \F_T$ as $\bbP_{t_i}^\gamma(A)$ with
\beal
\label{eqn27}
\bbP_{t_i}^\gamma(A) \doteq \bbP_0\bigg( &\w \in \O :\{ S_u(\w), t_i\leq u \leq T\} \in A;\\
&\textrm{ given }\{ S_r(\w):0\leq r \leq t_i\}\bigg).
\eal
We further take the convex hull of $\Qti$ and  denote it by $\Qti^c$. Namely, 
\beqn 
\label{eqn329}
\Qtic = \big\{ P \in \pft| P(A) = \sum_{i=1}^m \alpha_i P_i(A) \big\},
\eeqn 
for all $A \in \F_{T}$ with $P_i \in \Qti$, $0 \leq \alpha_i \leq 1$ with $\sum_{i=1}^m \alpha_i = 1$ for $i=1,\ldots, m$.

\subsection{Financial Scenario}
We consider the problem of investing in a risky asset $S_t$ and riskless asset $R_t = e^{rt}$, where $r > 0$ is the fixed interest rate. We assume that risky asset's dynamics as in \eqref{eqn250} is 
\beqn 
dS_t = \mu_t S_t dt + \sg_t S_t dW^S_t, \; \bbP_0\trm{-a.s.}
\eeqn 
For a given initial endowment $x_0 > 0$, the investor trades in a self financing way. Namely, denoting $\hat{\pi}_t$ as an $\F_t$ adapted stochastic process, which stands for the total amount of money invested in the risky asset $S_t$ at time $t$, $0 \leq t \leq T$, we have
\begin{align*}
dX^{\hat{\pi}}_t &= \hat{\pi}_t S_t^{-1} dS_t + (X^{\hat{\pi}}_t - {\hat{\pi}}_t)rdt \\
dX^{\hat{\pi}}_t &= {\hat{\pi}}_t (\mu_t dt + \sigma^S_t dW^S_t) + (X^{\hat{\pi}}_t - {\hat{\pi}}_t)rdt
\end{align*}

We further represent the amount of money invested in the risky asset as a fraction of current wealth via $\hat{\pi}_t =  X^\pi_t \pi_t$ for $0 \leq t \leq T$, where $\pi_t$ stands for the corresponding fraction at time $t$.  
Hence, for $X_0 = x_0$, the dynamics of wealth in this setting are given by 
\begin{align}
\label{eqn15}
d X_t^{\pi} &=  X_t^\pi \pi_t( \mu_t dt +  \sigma^S_t d W^S_t) + r X_t^\pi(1 - \pi_t) dt\\
X_T^\pi &= x_0\exp \int_0^T\{ \pi_s \mu_s + r(1 - \pi_s) - \frac{1}{2} \pis \}ds + \int_0^T \pi_s \sigma^S_s d W^S_s.
\end{align}

\subsection{Investor's Problem}
The investor is risk-averse and maximizes the minimal expected logarithmic utility over her set of r.c.p.m's representing alternative models until time horizon $T$, as follows. She is diffident about the underlying dynamics of the stock $S_t$ and wants to reevaluate her optimization problem at prespecifed fixed time epochs $0= t_0 < t_1 <\ldots <t_N= T$. At each time $t_i$ for $i = 0,\ldots, N-1$, we write the optimization problem of the investor as
\beqn 
\label{eqn2160}
V(t_i,x) = \supq \infq \big\{ \bextq \big[ \log(X_{T}^\pi) \big]  \big\},
\eeqn 
where $\bextq[\cdot]$ stands for the conditional expectation with respect to $Q \in \Qti$ as defined in \eqref{eqn26}. 
For time $t$ not equal to $t_i$'s, the investor sticks to her optimal $\gamma^{t_i}_{*} \in \Gamma$, assumes the model is correct until $t_{i+1}$ and solves the classical expected logarithmic utility problem based on
\beqn 
\label{eqn270}
V(t,x) = \supq \bbE^{\bbP_0}\big[ \log(X_{T}^\pi) | X^\pi_t =x \big]
\eeqn
Here, the process $S$ satisfies
\beqn 
dS_t = S_t(\mu_t dt) + \sg^S_t dW^S_t,
\eeqn 
where  
\beal
d\mu_t &= \theta^\mu_{t_i}(\eta^{\mu}_{t_i} - \mu_t)dt + \sg^{\mu}_{t_i}dW^\mu_t, \\
d(\sg^S_t)^2 &= \theta^\sg_{t_i}(\eta^\sg_{t_i}-(\sg^S_t)^2)dt + \xi(\sg^S_{*})^2dW^\sg_t.
\eal 
with $t_i < t < t_{i+1}$, where $\gamma^{t_i}_* = (\tmu_{t_i},\emu_{t_i},\tsg_{t_i},\esg_{t_i})$. Hence, the investor reevaluates her model by prespecified time intervals to be robust against the model risk.

Moreover, $\Pi_{\mathrm{ad}}$ stands for the set of admissible cash-values for the investor at time episode $\ti$, which is defined as follows.

\begin{definition}\label{def21}
Let $\pi = (\pi_s)_{\{ t \leq s \leq T\}}$ denote the $\B([t,T])\otimes \F_{T}$ progressively measurable process representing the cash-value allocated in the risky asset for $[t,T]$. We call $ (\pi_s)_{\{s \leq t \leq T\}}$ admissible and denote it by $\pi \in \Pi_{\mathrm{ad}}$, if
\beqn 
\bex\big[\int_{t}^{T} | \pi_s |^4  ds\big| \F_t] < \infty,\; \bbP_0\trm{-a.s.} 
\eeqn 
and $X^\pi_t > 0$ for all $t \leq s \leq T$, $\bbP_0$-a.s.
We say in this case that $\pi \in L^4([t,T])$ and note that $\Pi_{\mathrm{ad}}$ is nonempty, since $\pi \equiv 0$. We further say that $\pi_n \rar \pi$, if 
\beqn 
\bex\big[ \int_{t}^{T} | \pi_s - \pi^n_s | ^4 \big | \F_t]ds \rar 0,\; \bbP_0\trm{-a.s.},
\eeqn 
as $n \rar \infty$. We say in this case that $\pi_n$ converges to $\pi$ in $L^4([t,T])$.
\end{definition}
To attack problem \eqref{eqn2160}, we first solve \eqref{eqn270}. The optimization problem reads as
\beqn
\label{eqn124}
\supi \bextt[\log(X_T^{\pi})]. 
\eeqn
By Ito's formula, we get that
\begin{align} 
\label{eqn216}
&\supi \bextt[\log(X_T^{\pi})]  \\
&\indeq 
= \log(x_t) + \supit \bextt\bigg[ \int_t^{T} (\pi_s (\mu_s - r)+ r) 
\\ &\indeq\indeq 
-\frac{1}{2} \pis ds \bigg] 
\end{align}
By concavity on $\pi$, we conclude that checking first order condition inside the expectation on $\pi$ is sufficient and get that 
\begin{align}
(\mu_s - r ) - \sigma^2_s \pi_s &= 0.
\end{align}
Hence, we have 
\beqn
\label{eqn129}
\pi^*_s = \frac{\mu_s - r}{\sigma^2_s}
\eeqn
for all $t \leq s \leq T$, $\bbP_0$ a.s.,
and the optimal value function reads as
\beqn 
\label{eqn223}
\bextt[\log(X^\pi_{T})] = \log(x_t) + \bextt\big[\int_t^{T} r +\frac{(\mu_s -r)^2}{2\sigma^2_s}ds\big].
\eeqn

Going back to the robust optimization problem at each time $0 \leq t_i \leq T$ for $i = 0,\ldots,N$, we have
\beqn 
\label{eqn2320}
\sup_{\pi \in \Pi_{\mathrm{ad}}} \inf_{ Q \in \Q_{\ti}} \big\{ \bext\big[ \log(X_{T}^\pi) \big]  \big\}.
\eeqn 
To proceed, we first state our minimax result. 

\bethm \label{thm22}
Let $\Q^c_{\ti}$ be as in Equation \eqref{eqn329}, $\pi \in \Pi_{\mathrm{ad}}$ be as defined in Definition \ref{def21} and $X^\pi_t$ have the dynamics as in Equation \eqref{eqn250}. Then, we have 
\beqn 
\sup_{\pi \in \Pi_{\mathrm{ad}}}\inf_{Q \in \Q^c_{\ti}}\bigg\{\bext[\log(X^\pi_{T})] \bigg\}= \inf_{Q \in \Q^c_{\ti}}\sup_{\pi \in \Pi_{\mathrm{ad}}}\bigg\{\bext[\log(X^\pi_{T})]\bigg\}
\eeqn 
\ethm 
Theorem \ref{thm22} is an application of Sion's minimax theorem, which we recall here for convenience. 
\bethm \cite{S} \label{thm21} 
Let $X$ be a compact convex subset of a linear topological space and $Y$ a convex subset of a linear topological space. Let $f$ be a real-valued function on $X \times Y$ such that
\beitem
\item $f(x,\cdot)$ is upper semi continuous and quasi-concave on $Y$ for each $x \in X$. 
\item $f(\cdot, y)$ is lower semi continuous and quasi-convex on $X$ for each $y \in Y$.
\eitem 
Then 
\beqn 
\min_{x \in X}\sup_{y \in Y}f(x,y) = \sup_{y \in Y}\min_{x \in X}f(x,y). 
\eeqn 
\ethm 
We define first a suitable topology to work with the mapping 
\beqn 
P \rar \bbE^P_{t_i}[\log(X^\pi_{T})]\; \textrm{ for } P \in \Q^c_{\ti}. 
\eeqn 

\bedef \label{defn22}
A family of conditional probability measures  on $\O$, denoted by $\mathcal{S}_{\ti}$, is called relatively compact, if for every sequence $\{P_n\}$ in $\mathcal{S}_{\ti}$, there exists a subsequence $\{P_m\}$ of $\{P_n\}$ and a probability measure on $\O$ (not necessarily in $\mathcal{S}_{\ti}$) such that $\{P_m\}$ converges weakly to $P$. That is 
\beqn 
\bbE^{P_m}_{t_i}[g] \rar \bbE^{P}_{t_i}[g]
\eeqn 
for every $g:\O \rar \bbR$ with $g \in C_b(\O)$, where $C_b(\O)$ is the space of continuous bounded functions on $\O$. We say in this case that $P_m$ converges weakly to $P$, and denote it as $P_m \rhp P$. 
\edefn

\belem  \label{lem31}
$\Q^c_{\ti}$ as in Equation \eqref{eqn329} is compact with respect to topology defined in Definition \ref{defn22}. 
\elem 
\beproof
Let $\{P_n\}$ be a sequence in $\Q^c_{\ti}$. We will show that every sequence in  $\{P_n\}$ has a convergent subsequence in $\Q^c_{\ti}$. By taking a subsequence if necessary, one can assume that $P_n = \sum_{j=1}^m \alpha^j_n P^j_n$, where $\alpha^j_n \rar \alpha^j_{*}$  with $\alpha^j_{*} \in [0,1]$ and $P^j_n \in \Q_{\ti}$ for $j = 1,\ldots, m$. Recall that $P^j_n$ for $j=1, \dots, m$ are r.c.p.m's that are  induced by the processes of the form 
\beal
\label{eqn300}
&dS^j_T = \mu^j_TS^j_T dt + \sg^S_T S^j_T dW^j_T \; \bbP_0\;\textrm{a.s. for}\;s\leq t \leq T, \\
&\textrm{ given }\{S^{j}_r: 0\leq r \leq s\}
\eal
with $W^{S,j}_T$ is an $({\F_T}, \bbP_0)$ Wiener process as in Equation \eqref{eqn250}. Denoting $\varphi_{S^j_T}(\cdot)$ as the characteristic function of the conditional distribution induced by the process Equation \eqref{eqn300}, we have $P_n \rhp P^*$ in the sense of Definition \ref{defn22}, if and only if the characteristic functions $\varphi_{S_T}(z) =\bbE_{t}^{\bbP_0}[e^{izS_{T}^{j}}]$ converges pointwise to some characteristic function of some $\F_{T}$ measurable random variable $S$ for any $z \in \bbR$, where $S$ induces the probability measure $P^*$. 
But by dynamics of $(\sigma^2_t, \mu_t)$ as in Equation \eqref{eqnmu} and \eqref{eqnsg}, respectively, and since any sequence in $\Gamma$ has a convergent subsequence in $\Gamma$, this is only true if $(\theta_n^\mu,\theta_n^\sigma,\eta_n^\sg,\eta_n^\mu) \rar (\theta^\mu,\theta^\sigma,\eta^\sg,\eta^\mu)$. 

Hence, we have that the sequence of r.c.p.m.'s $P_n$ converges to some r.c.p.m. $P^*$ on $\O$, which is induced  by some $\F_{T}$ measurable random variable $S$, whose characteristic function is of the form
\beqn 
\sum_{j=1}^m \alpha^j_{*} \varphi_{S^j_T}(z) .
\eeqn 
Hence, we conclude the result. 
\eproof

We continue with the following lemma. 
\belem 
\label{cor31}Let $Q \in \Q^c_{\ti}$ with
\beqn
Q = \sum_{j=1}^m \alpha_j Q_j,
\eeqn
with $\alpha_j \geq 0$, $\sum_{j=1}^m \alpha_j = 1$ and $Q_j \in \Q_{\ti}$ for $j = 1,\ldots,m$, where $\Q_{\ti}$ is as in Equation \eqref{eqn26} 
and let $g: \O \rar \bbR$ be an $\F_{T}$ measurable mapping with $\bextq[|g(\w)|] < \infty$. Then, we have for any $m \in \bbN$
\beqn 
\bbE_{t_i}^{Q}[|g|] = \sum_{j=1}^m\alpha_j \bbE_{t_i}^{Q_j}[|g|]
\eeqn
with $\sum_{j=1}^m \alpha_j = 1$ and $\alpha_j \geq 0$. In particular, the mapping $Q \rar \bbE_{t_i}^Q[|\log(X^\pi_{T})|]$  is quasi-convex i.e. 
\beqn 
\bbE_{t_i}^{\alpha Q_1 + (1-\alpha)Q_2}[|\log(X^\pi_{T})|] \leq \max\{ \bbE_{t_i}^{Q_1}[|\log(X^\pi_{T})|], \bbE_{t_i}^{Q_2}[|\log(X^\pi_{T})|] \}
\eeqn 
\elem 
\beproof
Denoting the null sets of $Q_1,\ldots, Q_m$ as $N_{Q_1},\ldots, N_{Q_m}$ respectively, we note that $N_Q = \cap_{j=1}^m N_{Q_j}$. Without loss of generality,  by taking $g(\w) \geq 0$, $Q$-a.s. we get via an approximation of Lebesgue integration that
\beqn
\label{eqn332}
\sum_{j=1}^n c_j I_{A_j}(\w) \nearrow g(\w),\; Q\mathrm{-a.s.}
\eeqn 
for some real constants $c_j$, for $i =j,\ldots,n$, with $\cup_{j=1}^n A_j = \O$ and $A_j \cap A_k = \emptyset$ for $j \neq k$. Here, $I_{A_j}(\cdot)$ is the indicator function for the set $A_j \in \F_{t_i}$. Hence, we conclude via an approximation argument
\beqn 
\label{eqn333}
\bbE_{t_i}^Q[g] = \sum_{j=1}^m \alpha_j \bbE_{t_i}^{Q_j}[g]
\eeqn 
Since 
\beqn 
\max_{\gamma \in \Gamma}\bbE^{\bbP_0}_{t_i}[|\log(X^{\pi}_T)|] < \infty,
\eeqn 
we have $Q \rar \bbE^{Q}_{t_i}[\log(X^{\pi}_T)]$ is convex, in fact linear. Hence, we conclude the proof. 
\eproof

We continue with the following lemma. 
\belem \label{lem33}
Let $\pi$ and  $\pi_n$ be in $\Pi_{\mathrm{ad}}$ with $ \pi_n \rar \pi$ in $L^4(\O;[0,T])$ as $n \rar \infty$, $\pi_n(s) = \pi(s)$ for $0\leq s \leq t$ for all $n \in \bbN$, and let $P \in \Q^c_{\ti}$ be fixed. Then, the mapping
\beqn 
\pi \rar \bbE_{t_i}^P[\log(X^\pi_{T})]
\eeqn 
is continuous $\bbP_0$-a.s. Namely,
\beqn
\bbE_{t_i}^P[\log(X^{\pi_n}_{T})] \rar \bbE_{t_i}^P[\log(X^\pi_{T})] 
\eeqn
as $\pi_n \rar \pi$ in $L^4([t_i,T])$ $\bbP_0$ a.s. 
\elem

\beproof
By Lemma \ref{cor31}, we have that for any $P \in \Q^c_{\ti}$  
\beqn 
\bbE^P_{t_i}[\log(X^\pi_{T})] = \sum_{j=1}^m \alpha_j \bbE^{P_j}_{t_i}[\log(X^\pi_{T})]. 
\eeqn 
Hence, showing the continuity of 
\beqn 
\pi \rar \bbE^{P_j}_{t_i}[\log(X^\pi_{T})], 
\eeqn 
for $P_j \in \Q_{\ti}$ and hence for fixed $\gamma \in \Gamma$ is sufficient to show the continuity of $P \in \Q^c_{\ti}$. Since $\pi_n \rar \pi$ in $L^4([t_i,T])$, we note that 
\beqn
\sup_{n}\bext[\int_{t_i}^T (\pi^n_s)^4] < \infty,\;\bbP_0\;\mathrm{-a.s.}
\eeqn
and using Lemma \ref{lem21} we have 
\begin{align*}
&|\bext[\log(X^{\pi,\gamma}_{T})] - \bext[\log(X^{\pi_n,\gamma}_\piti)] |\\
&\indeq = |\bext[\int_{t}^\piti(\pi_s - \pi^n_s)(\mu_s - r) - \frac{1}{2}\bigg( \frac{\pi^2_s}{\sg^2_s} - \frac{(\pi^n_s)^2}{\sg^2_s} ds \bigg)]| \\
&\indeq \leq C\bext[\int_t^T \max_{\gamma \in \Gamma}(\mu_s - r)^2 ds]^{1/2}\bext[\int_t^T [\dpi^2ds]^{1/2}\\
&\indeq\indeq + C\bext[\int_t^T (\pi_s - \pi^n_s)(\pi_s + \pi^n_s) \frac{1}{\sg^2_s}ds]\\
&\indeq \leq \bext[\int_t^T (\pi_s+\pi^n_s)^4ds]^{1/2}\bext[\int_t^T \max_{\gamma \in \Gamma}\frac{1}{(\sg^S_s)^8}ds]^{1/2} + C\bext[\int_t^T (\pi - \pi^n_s)^2ds ] \\
&\indeq\rar 0, 
\end{align*}
as $n \rar \infty$. Thus, we have 
\beqn 
\bext[\log(X_{T}^{\pi_n,\gamma})] \rar \bext[\log(X_{T}^{\pi,\gamma})],
\eeqn
as $\pi_n \rar \pi$ for $n \rar \infty$ $\bbP_0$ a.s. Hence, we conclude the proof. 
\eproof

\belem
Let $\pi$ be in $\Pi_{\mathrm{ad}}$. Let $Q \in \Q^c_{\ti}$ be fixed. Then, the mapping
\beqn 
\pi \rar \bbE_{t_i}^Q[\log(X^\pi_{T})]
\eeqn 
is concave, in particular quasi-concave in $\pi$. 
\elem 
\beproof
By Lemma \ref{cor31}, it is enough to show the statement for $Q \in \Q_{\ti}$. Hence, for fixed $\gamma \in \Gamma$ by Ito lemma, we have 
\beqn 
\bext[\log(X_{T}^{\pi,\gamma})] =  \log(x_{t_i}) + \bextpi\bigg[ \int_{t_i}^{T} \pi_s(\mu_s -r) - \frac{1}{2}\pis ds \bigg], 
\eeqn 
from which it is easy to see that $\pi \rar \bext[\log(X_{T}^\pi)]$ is concave, hence quasi-concave in $\pi$. Hence, we conclude the result.  
\eproof 

Next, we continue with the following lemma. 
\belem \label{lem34}
Let $\pi \in \Pi_{\mathrm{ad}}$ be fixed and $Q,Q_n \in \Q^c_{\ti}$. Then, the mapping  
\beqn 
Q_n \rar \bbE_{t_i}^{Q_n}[\log(X^\pi_{T})]
\eeqn 
is  lower semi-continuous $\bbP_0$ a.s.. Namely,
\beqn 
\liminf_{n \rar \infty} \big\{\bbE_{t_i}^{Q_n}[\log(X^{\pi}_{T})]  \big\} \geq \bbE_{t_i}^Q[\log(X^\pi_{T})],
\eeqn as $Q_n \rhp Q$ in the sense of Definition \ref{defn22} as $n \rar \infty$, $\bbP_0$-a.s.
\elem 
\beproof
We have 
\beqn 
Q_n = \sum_{j=1}^m \alpha^j_n Q_n^j 
\eeqn 
where $0 \leq \alpha^j_n \leq 1$ and $Q^j_n \in \Q_{\ti}$ for $j=1,\ldots,m$. By Lemma \ref{cor31}, this means that 
\beqn
\bbE^{Q_n}_{t_i}[\log(X^\pi_{T})] = \sum_{j=1}^m \alpha^j_n \bbE^{\bbP_0}_{t_i}[\log(X^{\pi,\mu_n^j,\sigma^{S^j_n}}_{T})],
\eeqn
where $(\mu_n^j,\sigma_n^j)$ for $j = 1,\dots,m$ defines uniquely $P_j \in \Q^{\ti}$. 
Hence, we have 
\begin{align*}
\bbE^{P_n}_{t_i}[\log(X^\pi_{T})] &\leq \bbE^{P_n}_{t_i}[|\log(X^\pi_{T})|] \\
&\leq \bbE^{\bbP_0}_{t_i}[\max_{\gamma \in \Gamma}|\log(X^{\pi,\gamma}_{T})|]
< \infty,
\end{align*}
where 
\beqn 
\log (X^{\pi,\gamma}_{T}) = \log(x_{t_i}) + \big( \int_{t_i}^{T} \big(\pi_s(\mu_s - r) + r - \frac{1}{2}\frac{\pi_s^2}{\sg^2_s}ds \big) + \int_{t_i}^{T} \pis dW^S_s \big). 
\eeqn 
Next, we truncate our utility function $\log(x)$ as 
\beqn
V_k(x) =
\begin{cases} 
k &\textrm{if }\;\log(x) \geq k \\
\log(x) &\textrm{if }\;|\log(x)|\leq k \\
-k &\textrm{if }\;\log(x) \leq -k
\end{cases}
\eeqn
for $x > 0$ and for $k > 0$. 
We note that 
\beqn 
\bbE^{\bbP_0}_{t_i}[\log(X^{\pi,\gamma}_{T})] + \epsilon(k) \geq   \bbE^{\bbP_0}_{t_i}[V_k(X^{\pi,\gamma}_{T})],
\eeqn 
uniformly for all $\gamma \in \Gamma$ for some $\epsilon(k)$ depending on $k$ only with $\epsilon(k) \downarrow 0$ as $k \rar \infty$. Indeed, we have
\begin{align*}
\bbE_{t_i}^{P_0}[\log X^{\pi,\gamma}_{T} I_{\{|\log X^{\pi,\gamma}_{T}| > k\}}] \leq \max_{\gamma \in \Gamma } \bbE_{t_i}^{P_0}[\log X^{\pi,\gamma}_{T} I_{\{|\log X^{\pi,\gamma}_{T}| > k\}}] < \epsilon,
\end{align*}
for $k$ large enough. Furthermore, since by Lemma \ref{lem21}
\beqn 
\bbE^{\bbP_0}_{t}[\max_{\gamma \in \Gamma}|\log(X^{\pi,\mu,\sigma,r}_{T})|] < \infty\;\bbP_0\trm{-a.s.}
\eeqn 
is integrable with
\begin{align}
\label{eqn_Leb}
&V_k(x) \leq |\log(x)|\; \textrm{for any }k \geq 0\\
&\bbE_{t_i}^{\bbP_0}[V_k(X^\pi_{T})] \leq \bbE_{t_i}^{\bbP_0}[ \max_{\gamma \in \Gamma}|\log(X^\pi_{T})|] < \infty\; \bbP_0 \textrm{ a.s.},
\end{align}
we have 
\begin{align}
&\liminf_{n \rar \infty} \bbE^{P_n}_{t_i}[\log(X^\pi_{T})] = \liminf_{n\rar \infty} \bigg( \sum_{j=1}^m \alpha^j_n \bbE^{\bbP_0}_{t_i}[\log(X^{\pi,\mu^n_j,\sigma^n_j}_{T})]  \bigg) \\
&\quad \geq \liminf_{n\rar \infty} \bigg(  \sum_{j=1}^m \alpha^j_n\bbE^{\bbP_0}_{t}[V_k(X^{\pi,\mu^n_j,\sigma^n_j}_{T})] \bigg) - \epsilon(k)  \\
&\quad = \bbE^{P}_{t_i}[V_k(X^\pi_{T})] - \epsilon(k),\; \bbP_0 \textrm{ a.s.}
\end{align}
where the last equality is due to convergence $P_n \rhp P$. Finally, by letting $k \rar \infty$  with $\epsilon(k) \downarrow 0$ and via Lebesgue dominated convergence theorem by Equation \eqref{eqn_Leb}, we conclude the result. 
\eproof

\beproof[Proof of Theorem \ref{thm22}] We have that $\Pi_{\mathrm{ad}}$ is a convex set by Definition \ref{def21}. Similarly, $\Q^c_{\ti}$ as in Equation \eqref{eqn329} is a convex topological subset of $\cP(\O|\F_{t_i})$. Further, by Lemma \ref{cor31} we have that
\begin{align*}
\bbE_{t_i}^P[\log(X^\pi_{T})]
&= \sum_{j=1}^m \alpha_j \bbE^{\bbP_0}_{t_i}[\log(X^{\pi,\mu_j,\sigma_j}_{T})].
\end{align*}
Hence, $\bbE_{t_i}^P[\log(X^\pi_{T})]$ is real valued for any fixed $\pi \in \Pi_{\mathrm{ad}}$ and any $P \in \Q^c_{\ti}$. Further, by Lemma \ref{lem31}, \ref{lem33} and \ref{lem34}, we have that the conditions for Theorem \ref{thm21} are satisfied. Hence, we conclude the result via Sion's minimax theorem.
\eproof

\section{The Optimal Solution of the Investment Problem}
We have the following series of equations 
\begin{align}
\infit \supit \big\{ \bbE_{t_i}^{Q}[\log(X^{\pi}_{T})] \big\} &= \infir \supit \bbE_{t_i}^{Q}[\log(X^{\pi}_{T})] \\
&= \supit \infir \bbE_{t_i}^{Q}[\log(X^{\pi}_{T})]  \\
&\leq \supit \infit \bbE_{t_i}^{Q}[\log(X^{\pi}_{T})] \\
&\leq \infit \supit \bbE_{t_i}^{Q}[\log(X^{\pi}_{T})] 
\end{align}
where the first equality is by Lemma \ref{cor31}, the second equality is by Theorem \ref{thm22}, whereas the first inequality is by $\Q_{\ti}^c \supset \Q_{\ti} $ and second inequality is by classical minimax inequality.
Thus, we have 
\begin{align}
&  V (t_i,x,\nu,\sg)  = \infit \supit \big\{ \bbE^{Q}_{t_i}[\log(X^{\pi}_{T})]  \big\} \\
& V(t_i,\mu,\nu) = \log(x_{t_i}) + \inf_{\gamma \in \Gamma}  \{ \bbE^{Q}_{t_i} \big[ \int_{t_i}^{T} r + \frac{1}{2}\frac{(\mu_s - r)^2}{\sg^2_s} \big] ds.
\end{align}
For ease of notation supressing index $s$ for $\mu_{s}$ and $(\sg^S_s)^2$ and denoting $\nu = (\sg^S_s)^2$ the corresponding Hamilton-Jacobi-Bellman(HJB) equation reads as 
\beal
0 = \min_{\gamma \in \Gamma} \bigg( r &+  \frac{1}{2}\frac{(\mu-r)^2}{\nu} + V_t + V_\mu(\theta^\mu(\eta^\mu - \mu)) + V_\nu(\theta^\sg(\eta^\sg - \nu))\\
&+\frac{1}{2}V_{\mu\mu}(\sg_\mu)^2 + \frac{1}{2}V_{\nu\nu}\xi^2\nu^2
\bigg)
\eal 

\beitem
\item For $\mu < r$, hence $V_\mu < 0$
\benum
\item
If $\mu \in [\eta_{\min}^{\mu}, \eta_{\max}^{\mu}]$, choose $\eta^\mu = \eta^{\mu}_{\max}$ and $\theta^\mu = \theta^\mu_{\max}$.
\item 
If $\mu < \eta^{\mu}_{\min}$, choose $\eta^{\mu} = \eta^{\mu}_{\max}$ and $\theta^{\mu} = \theta^{\mu}_{\max}$.
\item 
If $\mu > \eta^{\mu}_{\max}$, choose $\eta^{\mu} = \eta^{\mu}_{\min}$ and $\theta^\mu = \theta^{\mu}_{\max}$. 
\enum
\item For $\mu > r$, hence $V_\mu > 0$ 
\benum 
\item
If $\mu \in [\eta_{\min}^{\mu}, \eta_{\max}^{\mu}]$, choose $\eta^\mu = \eta^{\mu}_{\min}$ and $\theta^\mu = \theta^\mu_{\max}$. 
\item 
If $\mu < \eta^{\mu}_{\min}$, choose $\eta^{\mu} = \eta^{\mu}_{\min}$ and $\theta^{\mu} = \theta^{\mu}_{\min}$. 
\item 
If $\mu > \eta^{\mu}_{\max}$, choose $\eta^{\mu} = \eta^{\mu}_{\min}$ and $\theta^\mu = \theta^{\mu}_{\max}$. 
\enum

\item
We have $V_\nu < 0$ for any $\nu > 0$. 
\benum
\item
If $\nu \in [\eta^{\sg}_{\min}, \eta^{\sg}_{\max}]$, choose $\eta^\sg = \eta_{\max}^\sg$ and $\theta^{\sg} = \theta^\sg_{\max}$. 
\item
If $\nu > \eta^{\sg}_{\max}$, choose $\eta^\sg = \eta^\sg_{\min}$ and $\theta^{\sg}= \theta^\sg_{\min}$. 
\item
If $\nu < \eta^\sg_{\min}$, choose $\eta^\sg = \eta^\sg_{\min}$ and $\theta^\sg = \theta^\sg_{\max}$.
\enum
\eitem 
Hence, we have that at each time interval $[t_i,t_{i+1})$, there exists a unique parameter set $\gamma^* = (\theta^\mu_*,\eta^\mu_*,\theta^\sg_*,\eta^\sg_*) \in \Gamma$ satisfying 
\beal
\label{fekat}
0 &=  r + \frac{1}{2}\frac{(\mu-r)^2}{\nu} + V_\mu(\theta^\mu(\eta^\mu - \mu)) + V_\nu(\theta^\sg(\eta^\sg - \nu))\\
&+\frac{1}{2}V_{\mu\mu}(\sg_\mu)^2 + \frac{1}{2}V_{\nu\nu}\xi^2\nu^2\\
\label{fekat2}
V(T,\mu,\nu) &= \log(x_T)
\eal
The solution to \eqref{fekat} and \eqref{fekat2} is well-known by Feynman-Kac 
\beal
V(t_i,\mu,\nu,x) = \log(x_{t_i}) + \bbE^M\bigg[\int_{t_i}^T r + \frac{1}{2}\frac{(\mu_s-r)^2}{\nu_s} ds\bigg|\mu_{t_i} = \mu,\nu_{t_i} = \nu \bigg]
\eal
subject to 
\beal
d\mu_t &= \theta^\mu_*(\eta^\mu_{*} - \mu_t)dt + \sg^\mu_*dW^M_t \\
d\nu_t &= \theta^\nu_*(\eta^\sg_{*}-\nu_t)dt + \xi(\sg^S_{*})^2dW^M_t,
\eal 
at $t = t_i$, where $W^M$ is a 2 dimensional Brownian motion under probability measure $M$, and the initial conditions are $\mu_{t_i} = \mu, \nu_{t_i} = \nu$ and $\bbE^M[\cdot]$ stands for the expectation taken with respect to $M$. Hence, the investor solves at each time $t_i$ using one of the end points of the corresponding intervals in $\Gamma$. We note that due to continuity of the state variables $\mu$ and $\nu$, the optimal parameter set $\gamma^*$ stays the same for some time interval, and then as the regions for the state variables change, there is a switch from one corner to the other. Hence, $(\gamma^*_t)_{0\leq t \leq T}$ being a piecewise constant function on $[0,T]$ belongs to $\Gamma$ by Lemma \ref{lem21}. Similarly, the optimal parameter $\pi^*_t$ reads as $\frac{\mu_t - r}{\sg^S_t}$, which is an element of $\Pi_{\mathrm{ad}}$ for all $0 \leq t \leq T$. Hence, the optimal controls $(\pi^*,\gamma^*)$ are as specifed in the admissible policy sets $(\Gamma, \Pi_{\mathrm{ad}})$ justifying the consistency of the optimal control problem \eqref{eqn2160}.

\section{Conclusion}
We have studied a logarithmic utility maximization problem, where there is uncertainty on drift $\mu$ and volatiliy $\sg$. We assume that the drift and volatility terms are OU and GARCH(1) processes, respectively. Hence, they take real values in noncompact subsets of $\bbR$. We assume that the corresponding parameters of the drift terms of $\mu$ and $\sg$ are uncertain but assumed to be estimated in some compact interval. We show that at each time $0 \leq t \leq T$, the optimal parameters are chosen on the right or left corners of the corresponding available parameter interval. Hence, we show that at each time, the optimal $\pi^*$ and the value function $V(t,x)$ are exactly as in the case, where there is no uncertainty on the parameters of the dynamics, but instead at each time $t$, the corresponding optimal parameters are given by corners of the available parameter interval.

\end{document}